\documentclass{amsart}

\usepackage{verbatim}

\usepackage{amsmath,amssymb,amsthm}
\usepackage{xy}
\usepackage{graphicx}
\numberwithin{equation}{section}
\newtheorem{theorem}[equation]{Theorem}
\newtheorem{proposition}[equation]{Proposition}
\newtheorem{corollary}[equation]{Corollary}
\newtheorem{lemma}[equation]{Lemma}
\newtheorem*{theorem*}{Theorem}

\theoremstyle{remark}

\newtheorem{remark}[equation]{Remark}

\newtheorem*{remark*}{Remark}

\def\Aut{\mathop{\operator@font Aut}}

\input xy
\xyoption{matrix}
\xyoption{arrow}

\begin{document}

\title[Boundedness for threefolds in
$\mathbb P^6$]{Boundedness for some rationally connected threefolds in
$\mathbb P^6$}
\author[Marian Aprodu]{Marian Aprodu}
\author[Matei Toma]{Matei Toma}

\address{{\it Marian Aprodu}: 
Romanian Academy, 
Institute of Mathematics "Simion Stoilow" 
P.O. Box 1-764, RO 014700, 
Bucharest, Romania}
\email{marian.aprodu@imar.ro}
\address{{\it Matei Toma}: Institut Elie Cartan, UMR 7502, Universit\'e de Lorraine, CNRS, INRIA, Boulevard des Aiguillettes, B.P. 70239, 54506 Vandoeuvre-l\`es-Nancy Cedex, France }
\email{Matei.Toma@univ-lorraine.fr}

\thanks{Marian Aprodu thanks IECN Nancy and IHES
Bures-sur-Yvette for hospitality during the preparation
of this work. Marian Aprodu was partly supported by LEA MathMode,
and by the CNCS-UEFISCDI grant PN-II-ID-PCE- 2011-3-0288, contract no. 132/05.10.2011. Matei Toma was partly supported
by LEA MathMode.}

\maketitle

\begin{abstract}
We prove boundedness of rationally-connected threefolds in $\mathbb P^6$
under some extra-assumptions.
\end{abstract}

\section{Introduction}

Hartshorne and Lichtenbaum conjectured that rational surfaces
in $\mathbb P^4$ form a bounded family. This famous
conjecture was proved in a more general setup by
Ellingsrud and Peskine \cite{Ellingsrud-Peskine}. 
Since then, the higher-dimensional case has become the
center of attention. 
Precisely, it is expected that there are finitely many
components of the Hilbert scheme parameterizing  smooth subvarieties $X$
{\em not} of general type in $\mathbb P^N$ of dimension $\mathrm{dim}(X)\ge N/2$
(since any variety $X$ can be embedded in $\mathbb P^{2\mathrm{dim}(X)+1}$
the problem of boundedness only makes sense in this range).
M. Schneider \cite{Schneider} has given an affirmative
answer, for $\mathrm{dim}(X)\ge (N+2)/2$.
Hence the remaining question is  whether or not the same is
true for $N\in\{2\mathrm{dim}(X)-1,2\mathrm{dim}(X)\}$. The
case of threefolds in $\mathbb P^5$ has been  settled in
\cite{BOSS}, which brings us to the next non-trivial case
of smooth threefolds $X$ in $\mathbb P^6$, the object of this note.
Recall that only finitely many components of the Hilbert scheme of smooth varieties can appear if the degree is bounded.

We prove:

\begin{theorem}
\label{thm:main}
Let $X$ be a smooth rationally connected
threefold in $\mathbb P^6$ not contained in any
fourfold of degree $\le 34$ and denote by $S$ a general hyperplane section. If $K_S^2\le 9$, then 
the degree of $X$ is bounded by $34^3$.
\end{theorem}

The condition $K_S^2\le 9$ is realised if either $X$ is covered by lines, or $|K_X+H|=\emptyset$, where $H$ denotes the hyperplane class, or $S$ is not of general type, Corollary \ref{cor:main}. The hypothesis on rational connectedness will be used only through its consequences on cohomology vanishing, i.e. $h^i(\mathcal O_X)=0$ for $i\ge 1$.

Some other boundedness results in the same spirit
have been obtained by Sabatino \cite{Sabatino}. Precisely, in 
loc. cit. it is assumed that the general hyperplane
section of $X$ is a ruled surface.

It follows directly from \cite{Ionescu-Toma} that {\em subcanonical} varieties
not of general type form a bounded family. Hence we are mainly
concerned with the case of non-subcanonical subvarieties, although this
hypothesis will not be used.

The main technical ingredients used in the proof are the following: 
the lifting theorem of Chiantini-Ciliberto \cite{Chiantini-Ciliberto93}, 
the bounds on the genus of curves in $\mathbb P^4$ of Chiantini-Ciliberto-di Gennaro \cite{Chiantini-Ciliberto-DiGennaro}, 
the semi-positivity of Schur polynomials of Fulton-Lazarsfeld \cite{Fulton-Lazarsfeld},
and the Hodge index formula for some divisors in $X$.

\section{Notation,  setup}

The general setup is the following. 

\medskip

$X\subset \mathbb P^6$ is a rationally connected smooth variety of dimension 3,

$H:=\mathcal O_X(1)$ the hyperplane section bundle on $X$,

$K:=K_X$ the canonical bundle of $X$,

$h:=c_1(H)\in H^2(X,\mathbb Z)$ the class of $H$,

$k:=c_1(K)\in H^2(X,\mathbb Z)$ the class of $K$,

$N:=N|_{X|\mathbb P^6}$ the normal bundle of $X$,

$S\in |H|$ a general hyperplane section,

$H_S:=H|_S$ the induced very ample bundle on $S$,

$C\in |H_S|$ a general sectional curve,

$H_C:=H|_C$ the induced very ample bundle on $C$,

$d:=(H^3)$ the degree of $X$,

$g:=g(C)$ the sectional genus,

$\delta:=2g-2=h^2\cdot k+2h^3$ the degree of $K_C$,

$c_i:=c_i(X)\in H^{2i}(X,\mathbb Z)$, $i\in\{1,2,3\}$ the Chern classes of $X$,

$n_i:=c_i(N)\in H^{2i}(X,\mathbb Z)$, $i\in\{1,2,3\}$ the Chern classes of the normal bundle,

$\chi:=\chi(\mathcal O_S)$ the Euler characteristic of $S$,

$u:=h^{1,1}(S)$ the Picard number of $S$,

$v:=c_1(N(-H))^3$.

\medskip

Note that $h^i(\mathcal O_X)=0$ for all $i\ge 1$ \cite{Kollar}, Corollary IV.3.8.

\section{Identities}

In this section we record a number of useful relations in connection with our setup.

\begin{lemma}
\label{lemma:q=0}
 $S$ is a regular surface i.e. $q(S)=0$. Moreover
$\chi=\chi(K+H)+1=h^0(K+H)+1$
and hence $p_g(S)=0$ if and only if $|K_X+H|=\emptyset$.
\end{lemma}

\proof
 Apply $h^2(\mathcal O_X(-1))=0$ (from Kodaira vanishing) 
and the vanishing of $h^1(\mathcal O_X)$ and $h^2(\mathcal O_X)$
to the exact sequence
\[
 0\to \mathcal O_X(-1)\to \mathcal O_X\to \mathcal O_S\to 0.
\]
\endproof

\begin{lemma}
 \label{lemma:K_S^2} $c_2(S)=2\chi+u$ and
$K_S^2=10\chi-u$.
\end{lemma}

\proof
From the Gauss-Bonnet formula, we have
$c_2(S)=\chi_{\mathrm{top}}(S)=2-2b_1(S)+b_2(S)$,
from where, applying Lemma \ref{lemma:q=0},
we obtain $c_2(S)=2+b_2(S)=2\chi+u$.

For the second relation, apply Noether's formula on $S$.
\[
 \chi=\frac{1}{12}(K_S^2+\chi_{\mathrm{top}}(S)),
\]
which implies $K_S^2=12\chi-(2\chi+u)$.
\endproof

\begin{lemma}
 \label{lemma:kc_2}
$k\cdot c_2=-24$.
\end{lemma}

\proof
From the Hirzebruch-Riemann-Roch formula, we have
$\chi(\mathcal O_X)=c_1\cdot c_2/24$, and from rational-connectedness
we have $\chi(\mathcal O_X)=1$.
\endproof

\begin{lemma}
 \label{lemma:c(N)}
The Chern classes of $N$ are the following:
\begin{enumerate}
 \item $n_1=7h+k$;
 \item $n_2=21h^2+7h\cdot k+k^2-c_2$;
 \item $n_3=35h^3+21h^2\cdot k+7h\cdot k^2+k^3-7h\cdot c_2-c_3+48$.
\end{enumerate}
\end{lemma}

\proof
We use the normal sequence
\begin{equation}
 \label{eqn:normal}
0\to T_X\to T_{\mathbb P^6}|_X\to N\to 0
\end{equation}
which implies that $c(N)=c(T_{\mathbb P^6}|_X)/c(X)$.
Obviously
\[
 \frac{1}{c(X)}=1-c_1t+(c_1^2-c_2)t^2+(-c_1^3+2c_1c_2-c_3)t^3
\]
i.e.
\[
 \frac{1}{c(X)}=1+kt+(k^2-c_2)t^2+\left(k^3+48-c_3\right)t^3
\]
On the other hand, since $c(T_{\mathbb P^6})=(1+ht)^7$, we obtain
\[
 c_1(T_{\mathbb P^6}|_X)=7h,\ c_2(T_{\mathbb P^6}|_X)=\binom{7}{2}h^2=21h^2,\
 c_3(T_{\mathbb P^6}|_X)=\binom{7}{3}h^3=35h^3.
\]
This implies
\[
 c(N)=(1+7ht+21h^2t^2+35h^3t^3)\cdot 
\left(1+kt+(k^2-c_2)t^2+\left(k^3+48-c_3\right)t^3\right).
\]
\endproof

\begin{lemma}
\label{lemma:dp}
 $n_3=d^2$.
\end{lemma}

\proof
This is the double-point formula \cite{Lascu-Mumford-Scott}.
\endproof

\begin{lemma}
\label{lemma:rel}
 We have the following numerical relations:
\begin{enumerate}
 \item $h^2\cdot k=-2d+\delta$;
 \item $h\cdot k^2=3d-2\delta +10\chi-u$;
 \item $k^3=-4d-24\delta-120\chi+12u+v$;
 \item $h\cdot c_2=d-\delta+2\chi+u$.
 \item $c_3=3d-10\delta-64\chi-2u+v-d^2+48$;
\end{enumerate}
\end{lemma}

\proof

For (1), by adjunction we have $\delta=(k+2h)\cdot h^2$,
and hence $h^2\cdot k=-2d+\delta$.

For (2) we  apply Lemma \ref{lemma:K_S^2}:
$K_S^2=10\chi-u$.
On the other hand $$K_S^2=(k+h)^2\cdot h=h^3+2h^2\cdot k+h\cdot k^2$$
implying $h\cdot k^2=10\chi-u-d-2(-2d+\delta)$.

For (3) we use the definition of $v$. Note that
$c_1(N(-1))=4h+k$, from Lemma \ref{lemma:c(N)},
and hence
\[
 v=64h^3+48h^2\cdot k+12h\cdot k^2+k^3.
\]
The preceding formulae show that
\[
 k^3=v-64d-48(-2d+\delta)-12(3d-2\delta +10\chi-u).
\]

For (4) remark that $h\cdot c_2=c_2(T_X|_S)$.
The exact sequence
\[
 0\to T_S\to T_X|_S\to N_{S|X}\to 0
\]
implies that (note that $N_{S|X}=H_S$) 
$h\cdot c_2=c_2(S)+c_1(S)\cdot H_S$.
From Lemma \ref{lemma:K_S^2} it follows that
$h\cdot c_2=2\chi+u-(k+h)\cdot h^2$.

For (5), apply Lemma \ref{lemma:c(N)}, Lemma \ref{lemma:dp}
and the relations above:
$c_3=35h^3+21h^2\cdot k+7h\cdot k^2+k^3-7h\cdot c_2-n_3+48$
$=35d+21(-2d+\delta)+7(3d-2\delta +10\chi-u)+(-4d-24\delta-120\chi+12u+v)
-7(d-\delta+2\chi+u)-d^2+48$.
\endproof

\section{Inequalities}

This section is devoted to some useful inequalities
between the given invariants.

\subsection{A Lifting Theorem}\cite{Chiantini-Ciliberto93}
It is known \cite{Roth} that if the sectional curve $C$ is contained in 
a hypersurface of degree $s$ 
in $\mathbb P^4$, with $s^2<d$, then both $S$ and $X$ are contained in
hypersurfaces of degree $s$ in $\mathbb P^5$, and $\mathbb P^6$
respectively.

We shall need a version of this result, 
which is a special case of Theorem 0.2 \cite{Chiantini-Ciliberto93}.

\begin{theorem}
\label{thm:lifting}
 If the general sectional curve $C$ of $X$ is contained in a surface
 in $\mathbb P^4$ of degree $s$ such that $d>\frac{(s-1)(s-3)}{2}+8s-3$ 
 then $X$ is contained in a 4-fold of degree $s$ in $\mathbb P^6$.
\end{theorem}

This result will be used together with the bounds of the 
genus of $C$, see below.

\subsection{Bounds on the genus of $C$}
We shall apply the main result of \cite{Chiantini-Ciliberto-DiGennaro}
in our setup; the general sectional curve $C$ is non-degenerate.

\begin{theorem}
\label{thm:CCdG}
Notation as above. Assume that $C$ is not contained in a 
surface of even degree $s\ge 11$ with $d>s^3$. Then the genus
of $C$ is bounded by
\[
\frac{d^2}{2s}+\frac{d}{2}\left(\frac{s}{2}-3\right)+\frac{3s^2-20}{8}.
\]
\end{theorem}

When making the substitutions in the main result of \cite{Chiantini-Ciliberto-DiGennaro}
we do not take care of optimality of the bound.
A similar bound can be found for odd $s$, either
by applying directly loc. cit. or replacing $s$ by $s-1$
in the statement above.

\subsection{Schur polynomials}
The bundle $N(-1)$ is globally generated, hence 
by \cite{Fulton-Lazarsfeld}, the associated Schur
polynomials 
\[
 s_{(1)}:=c_1(N(-1)),\ s_{(20)}:=c_2(N(-1)),\ s_{(300)}:=c_3(N(-1)),\
\]
\[
s_{(11)}:=(c_1^2-c_2)(N(-1)),\ s_{(210)}:=(c_1\cdot c_2-c_3)(N(-1))
\]
\[
s_{(111)}:=(c_1^3-2c_1\cdot c_2+c_3)(N(-1))
\]
are semi-positive, in particular $s_{(1)}\cdot h^2$,
$s_{(20)}\cdot h$ and $s_{(11)}\cdot h$ are non-negative.
We compute all these non-negative numbers.

\begin{lemma} Notation as above. We have
 \begin{enumerate}
  \item $s_{(1)}\cdot h^2=2d+\delta$;
  \item $s_{(20)}\cdot h=2d+4\delta+8\chi-2u$;
  \item $s_{(11)}\cdot h=d+2\delta +2\chi+u$;
  \item $s_{(300)}=-5d-5\delta-8\chi+2u+d^2$;
  \item $s_{(210)}=4d-3\delta-30\chi-3u+v+24-d^2$;
  \item $s_{(111)}=-3d+11\delta+68\chi+4u-v-48+d^2$.
 \end{enumerate}
\end{lemma}

\proof
We easily show that
$c_1(N(-1))=n_1-3h$, $c_2(N(-1))=n_2+3h^2-2h\cdot n_1$,
$c_3(N(-1))=n_3-h^3+h^2\cdot n_1-h\cdot n_2$.
We apply next Lemma \ref{lemma:c(N)}, \ref{lemma:kc_2},
\ref{lemma:dp}, \ref{lemma:rel}.
\endproof

\subsection{Hodge-index Theorem}

\begin{proposition}
\label{prop:Hodge}
 Let $D$ be a smooth divisor on $X$.
Then $$(H\cdot K\cdot D)^2\ge(H^2\cdot D)(K^2\cdot D).$$
\end{proposition}

\proof
Write $(H\cdot K\cdot D)=(H|_D\cdot K|_D)$, 
$(H^2\cdot D)=(H|_D^2)$, $(K^2\cdot D)=(K|_D^2)$
and apply the usual Hodge-index theorem on $D$ to $H|_D$ and $K|_D$.
\endproof

\begin{remark}
 Applying Hodge theorem for a linear combination $aD_1+bD_2$
of divisors does not give anything new apart
from the inequalities obtained for $D_1$ and~$D_2$.
\end{remark}

We can apply Proposition \ref{prop:Hodge} to two different
divisors on $X$. One is $H$ itself and the other one is $4H+K$;
recall that $\mathrm{det}(N(-1))=4H+K$ is globally generated.

\begin{corollary}
\label{cor:Hodge}
 $(3d+6\delta+10\chi-u)^2\ge v(2d+\delta)$ and
$\delta^2-(2\delta + 10\chi -u)d+d^2\ge 0$.
\end{corollary}

\section{Boundedness of the degree}

\noindent
{\bf Proof of Theorem \ref{thm:main}.}

Suppose that $d>34^3$.
From the hypothesis, we know that 
$10\chi-u\le 9$. Note that 
\[
s_{(20)}\cdot h+s_{(11)}\cdot h=3d+6\delta+10\chi-u\ge 0.
\]
From corollary \ref{cor:Hodge}, we
obtain
\[
v\le \frac{(3d+6\delta+10\chi-u)^2}{2d+\delta}\le \frac{(3d+6\delta+9)^2}{2d+\delta}.
\]
On the other hand, from the non-negativity of $s_{(210)}$, we
have
\[
v\ge d^2-4d+3\delta+30\chi+3u-24\ge d^2-4d+3\delta+9
\]
whence
\[
33\delta^2+(-d^2+34d+99)\delta+(81+17d^2+36d-2d^3)\ge 0
\]
Since for $d> 34^3$, the expression 
$81+17d^2+36d-2d^3$ is clearly negative, we have
\[
\delta\ge \frac{d^2-34d}{33}-3.
\]

From theorem \ref{thm:lifting} and
theorem \ref{thm:CCdG} applied for $s=34$, it follows that
\[
\delta=2g-2\le \frac{d^2}{s}+d\left(\frac{s}{2}-3\right)+\frac{3s^2}{4} -7=
 \frac{d^2}{34}+14d+860.
\]
These two opposite inequalities for $\delta$ yield
\[
 \frac{d^2}{34}+14d+860\ge \frac{d^2-34d}{33}-3
\]
which contradicts the assumption $d>34^3$. 
\qed

\begin{corollary}
\label{cor:main}
Let $X$ be a smooth rationally connected
threefold in $\mathbb P^6$ not contained in any
fourfold of degree $\le 34$. Then 
the degree of $X$ is bounded by $34^3$ if one of the
following conditions is satisfied:
\begin{enumerate}
 \item $X$ is covered by lines;
 \item the general hyperplane section $S$ of $X$ is not of general type;
 \item $|K_X+H|=\emptyset$, where $H$ denotes the hyperplane class.
\end{enumerate}
\end{corollary}

\proof
We prove that in each of the three cases, we have $K_S^2\le 9$.

If $S$ is not of general type, then the conclusion follows from
the classification of surfaces (note that $K^2$ decreases when we blow up).

If $X$ is covered by lines, then $|K_X+H|$ is empty. Indeed, 
by \cite{Ionescu-Russo}, it follows that for $\ell$ a general line
in a covering family, the restriction ${\Omega^1_X}|_\ell
\cong \mathcal O(-2)\oplus \mathcal O(a)\oplus \mathcal O(b)$
with $-1\le a,b\le 0$ (i.e. lines are free, standard). Hence
$K_X\cdot \ell\le -2$ and, since $H\cdot \ell=1$, $(K_X+H)\cdot \ell\le 1$.
Hence $K_X+H$ cannot be effective.

Under the assumption $|K_X+H|=\emptyset$, from Lemma \ref{lemma:q=0},
we infer that $p_g(S)=0$, and hence $\chi=1$, which implies that $10\chi-u\le 9$.
\endproof

\end{document}